\documentclass{tMAA2e}
\usepackage[utf8]{inputenc}
\usepackage{pgf,tikz}
\usetikzlibrary{calc}
\newcommand{\radius}{1.5cm}
\begin{document}
\title{Textile D-forms and $D_{4d}$}

\author{Katherine A.~Seaton$^{a}$$^{\ast}$\thanks{$^\ast$ Email: k.seaton@latrobe.edu.au
\vspace{6pt}} \\\vspace{6pt} $^{a}${\em{Dept. Mathematics and Statistics, La Trobe University VIC 3086, Australia}}
 }

\maketitle

\begin{abstract}
D-forms have in the past been created from inflexible materials, or considered as abstract mathematical objects. This paper describes a number of realisations of D-forms, and the related pita-forms, in textiles. Examples are given in which the created 3D space is purposefully employed and others in which ornamentation of the constituent surfaces is the highlighted feature. In particular, a set of biscornu has been fashioned to provide a 3D sampler of the axial point group $D_{4d}$ and its subgroups, using hitomezashi.


\begin{classcode} 20D05; 52A15; 97M80 \end{classcode}

\end{abstract}

\section{Introduction}
A \textit{D-form} is a 3-dimensional domain created when two isoperimetric developable sheets are joined together by seaming their boundary curves. The same two constituent sheets can give rise to many possible D-forms, by changing the point on each boundary at which the seaming process begins \cite{Wi}. First dreamt into concrete form by the designer Tony Wills  (literally dreamt, and literally concrete \cite{WW}) some 20 years ago, the idea has subsequently been placed on a firmer mathematical footing \cite{DO, DP, PW, Sh1}.

 The mathematical results, chiefly that a D-form is the convex hull of its seam curve, and that they have no creases (apart from the seam itself), require the refining condition that the constituent surfaces each have a smooth and convex  boundary curve \cite{DO}. (That is, the examples of D-forms in the design literature include objects to which these results do not apply.) Two generalisations of  D-forms proposed by mathematicians are \textit{pita-forms} (formed from 
 one developable sheet only, seamed to its own smooth convex boundary) \cite{DO} and \textit{seam-forms} (more than two such sheets are used) \cite{DP}. 

Of course, mathematicians work with ideal surfaces and curves; the material from which their D-forms are created does not stretch or collapse and their seams do not pucker. In the real world, stiff paper, card  or plastic (like overhead projector transparencies) work well as the developable sheets \cite{LB, KST, Sh2}.  Wills has made D-forms from plywood or metal both as artistic forms and in professional practice for street architecture \cite{Wi,WW}. To design a large final product from scratch, computational tools can be used to match edge length \cite{OWBSKBA}, and to predict the outcome \cite{PB, GAS}, which is otherwise generally hard to anticipate. There is, however, much delight in observing how the 3D object is formed progressively as the edges are brought together. An oft-quoted statement by Wills is to the effect that there is no ugly D-form \cite{Sh2}. In a recent paper, the authors ask whether we perhaps find this beauty because the human senses (unconsciously) appreciate their minimal convex nature \cite{BMR}.

There is no imperative (apart from the mathematical) to seam together only shapes with smooth convex boundaries, and indeed the earliest D-forms used squares and shapes with undulating edges \cite{Sh1, Sh2,Wi}. Recently, specified hierarchies from within the sphericon family of 3D objects were shown to be D-forms or pita-forms formed by seaming non-convex shapes \cite{KS}. The knitted and crocheted surfaces of the objects illustrating that paper were stabilised by felting or lining, or were made from a stiff plastic thread and lightly filled. Both Sharp \cite{Sh1} and Wills \cite{Wi} advise D-forms {\em not} be made of fabric, as it may shear or deform, though concede that this might be an interesting line to pursue.  In fact, in a commonplace context which predates Wills' dream by at least a century, it is precisely the deformation of fabric over a spherical core which creates a baseball from two non-convex leather `flats', joined with curved seam that pitchers use to their advantage \cite{BB}.

In this paper, some other extant textile D-forms are discussed. In the first section, designs for two pita-form accessories (a hat and a handbag) are described; the feature used purposefully in both is the internal 3D space. In the second section, the D-form cushions called \textit{biscornu} are introduced, and the customary decoration (embroidery) of their surfaces is used to illustrate their underlying symmetry group, that of the square anti-prism\footnote{There is a another set of objects called 15-sided biscornu, made from 15 squares which can be joined together in 3 different configurations; we are not discussing these.}. Strikingly, the hat design and the emergence of biscornu can both be dated to the first decade after Wills proposed D-forms, though there seems to be no evidence that this is not purely coincidence.

\section{Pita-form accessories}
\subsection{A pita-form handbag} 
\begin{figure}
\begin{center}
{\resizebox*{10cm}{!}{\includegraphics{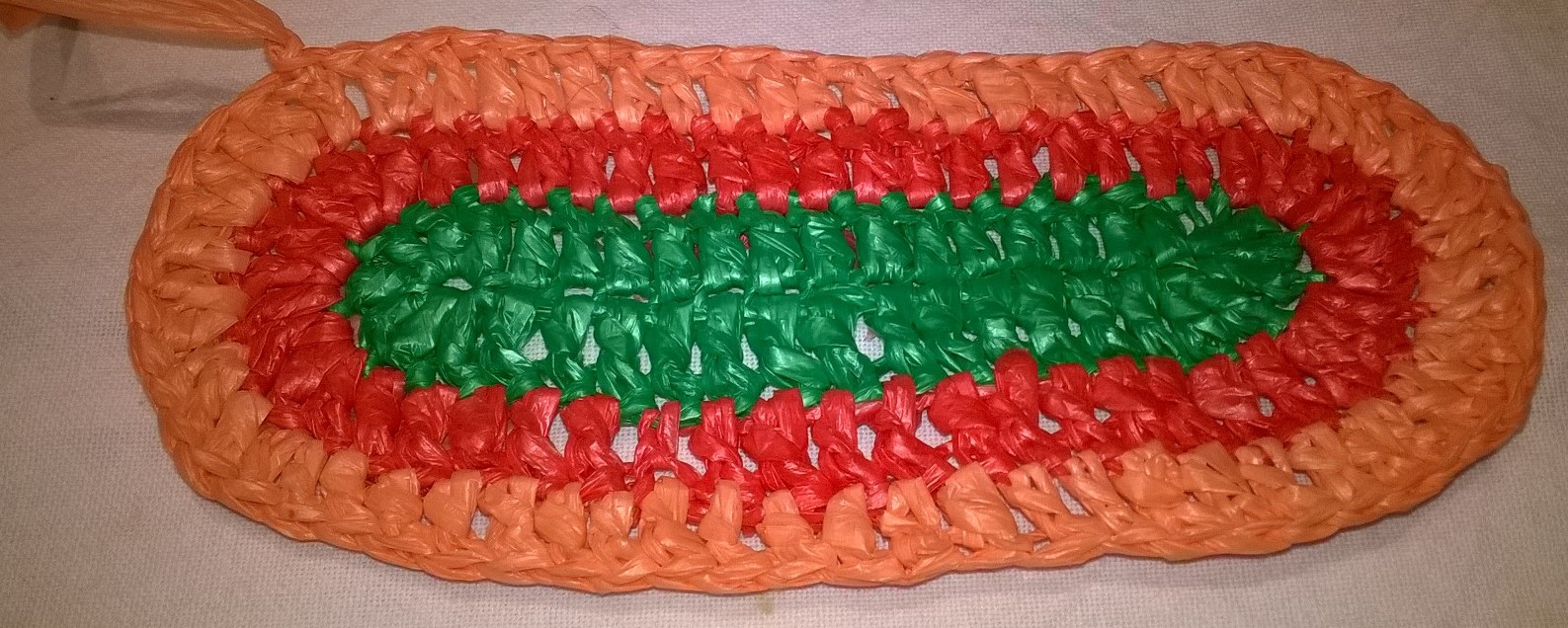}}}
\caption{\label{stadium} A crocheted stadium, executed in plarn (yarn made from plastic bags) by the author in 2015.}
\end{center}
\end{figure}

\begin{figure}
\begin{center}
\resizebox{5 cm}{!}{
\begin{tikzpicture}
  \coordinate (A) at (0,0);
  \coordinate (B) at (0,2*\radius);
  \coordinate (C) at (2*\radius,2*\radius);
  \coordinate (D) at (2*\radius,4*\radius);
  \draw (A) -- node[left] {$2r$} (B);
  \draw (B) -- node[above] {$2r$} (C);
  \draw (C) -- node[right] {$2r$} (D);
  \draw[densely dashed] (A) -- (C);
  \draw[densely dashed] (B) -- (D);
  \draw (D) arc [start angle=90, end angle=180, radius=2*\radius]
    arc [start angle=90, end angle=270, radius=\radius]
    arc [start angle=270, end angle=360, radius=2*\radius]
    arc [start angle=-90, end angle=90, radius=\radius];
  \draw ($(B)+(0,-0.2*\radius)$) -- ++(0.2*\radius,0) -- ++(0,0.2*\radius);
  \draw ($(C)+(-0.2*\radius,0)$) -- ++(0,0.2*\radius) -- ++(0.2*\radius,0);
\end{tikzpicture}
}
\caption{\label{trisphericon} The non-convex from which a trisphericon can be obtained as a pita-form. This figure also appears in \cite{KS}.}
\end{center}
\end{figure}
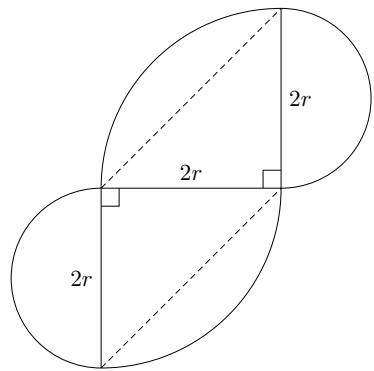

The investigation of sphericons recounted in \cite{KS} was sparked by a crocheted coin purse,  a pita-form obtained from a \textit{stadium} (see Figure \ref{stadium}). 
In Figure \ref{trisphericon}, the non-convex shape from which the \textit{trisphericon} is formed is shown for comparison. Its  curved edge is formed from quarter- and semi-circles.

\begin{figure}
\begin{center}
{\resizebox*{10cm}{!}{\includegraphics{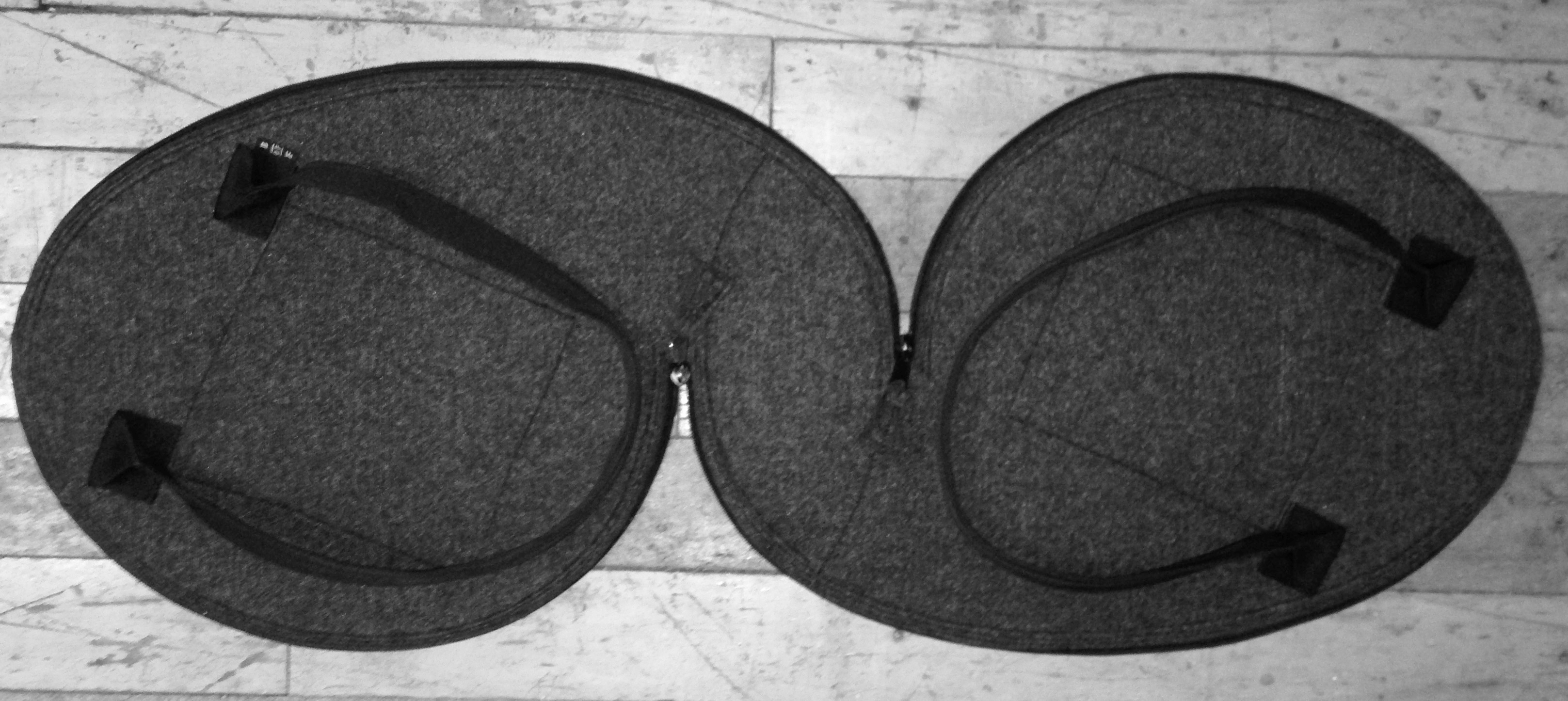}}}
\caption{\label{flat} The flat form devised by Atkinson for a pita-form handbag \textit{circa} 2016. Courtesy Bruce Atkinson.}
\end{center}
\end{figure}

\begin{figure}
\begin{center}
{\resizebox*{6cm}{!}
{\includegraphics{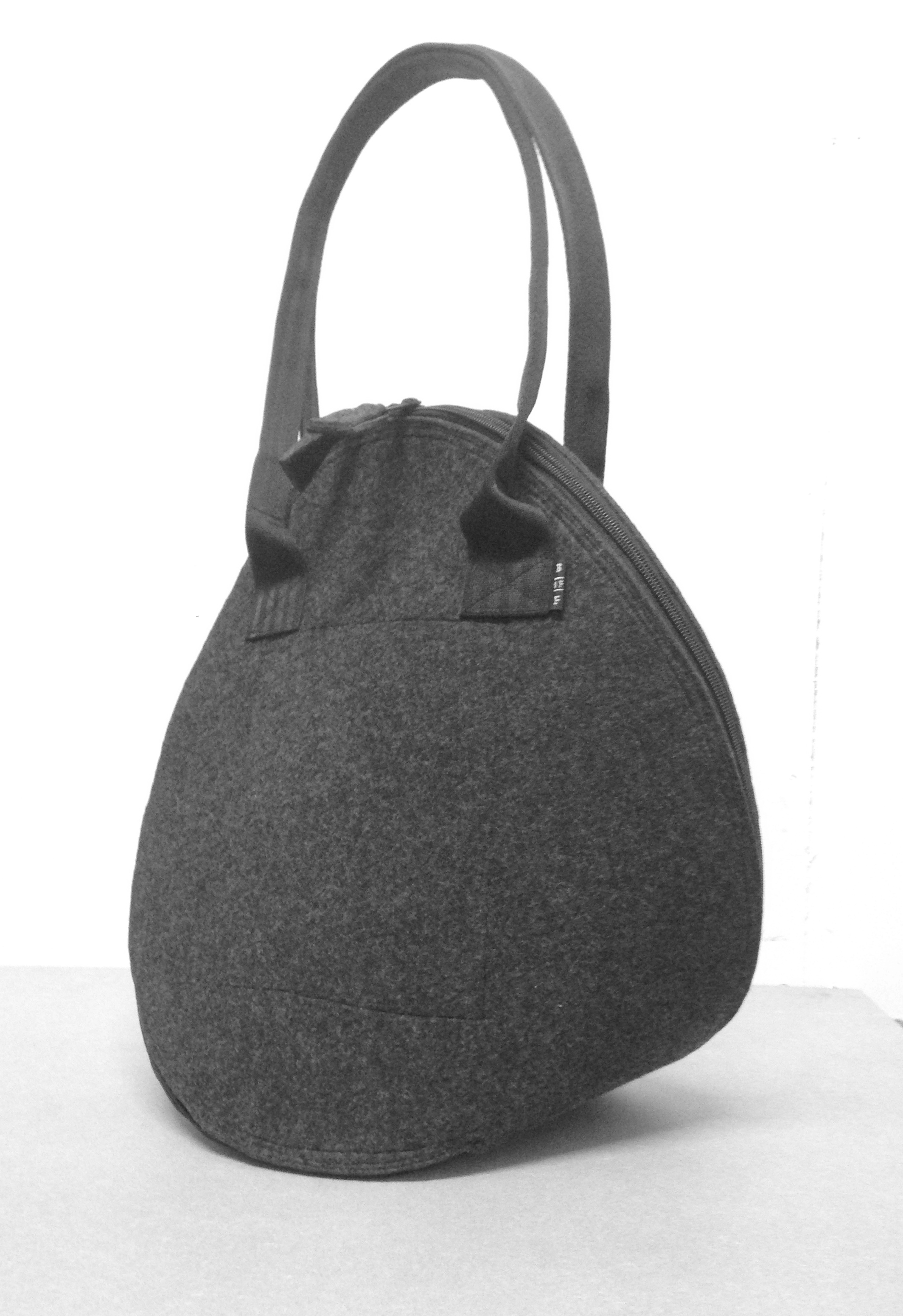}}
}
\caption{\label{up} The pita-form bag itself; is is large enough to hold an A4 notepad. Courtesy Bruce Atkinson.}
\end{center}
\end{figure}

As a result of \cite{KS}, an interesting correspondence ensued with UK-based pattern-cutter and fashion designer Bruce Atkinson. Investigating new possible garment forms \cite{BA}, he had been experimenting by D-forming (though he did not know it under that name at the start of our correspondence). He has given permission to include here images of an unusual handbag of his design, which can be seen both in its developed and seamed (zipped) form in Figures \ref{flat} and \ref{up}. Looking at Figure \ref{flat}, it is apparent why Figure \ref{trisphericon} caught his attention. The curved edges in Atkinson's work are found by informed trial-and-error \cite{BA}.

\subsection{Almost a pita-form hat} 
\begin{figure}
\begin{center}
{\resizebox*{8cm}{!}
{\includegraphics{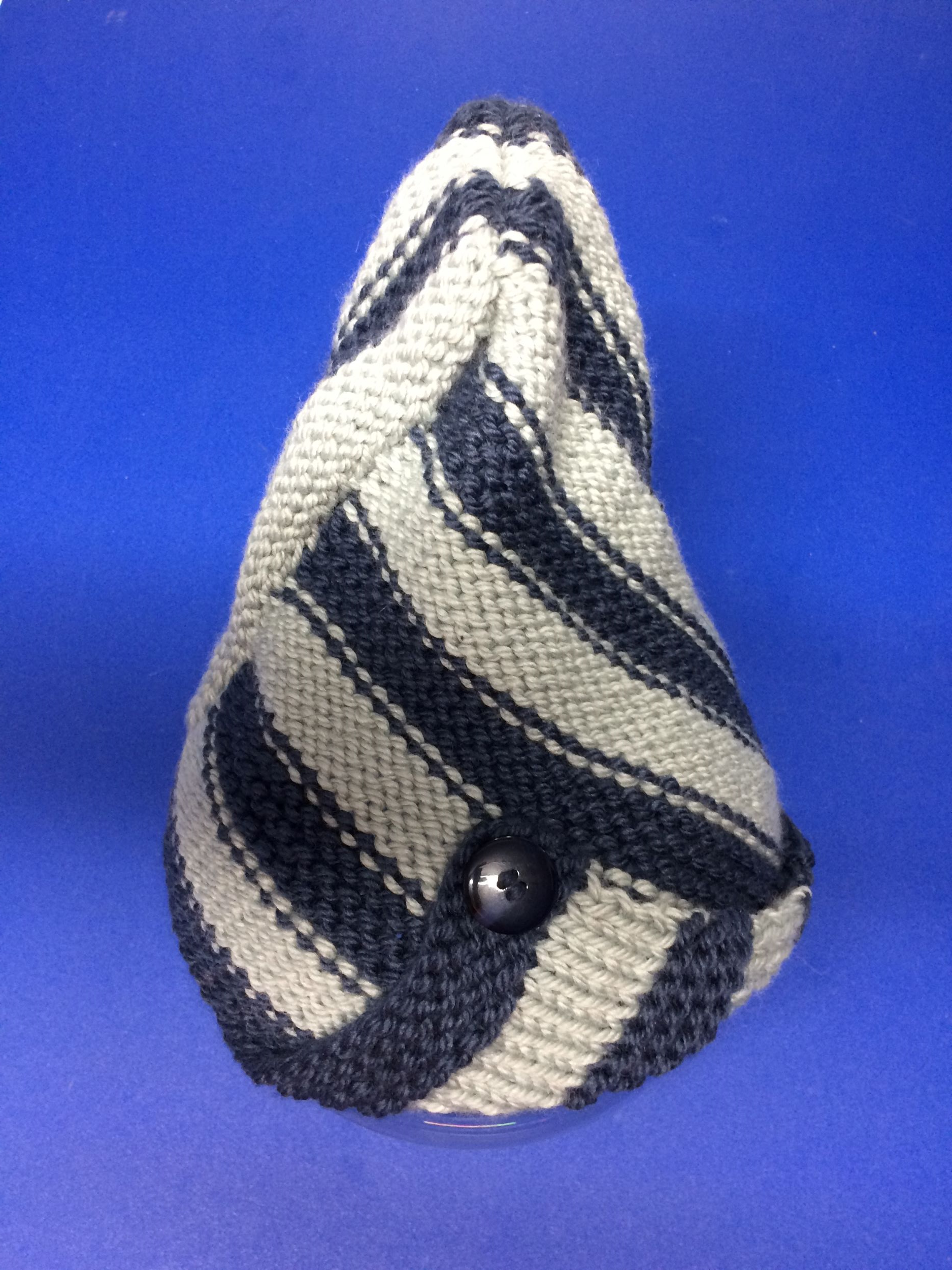}}
}
\caption{\label{orpheus} \textit{Orpheus} hat, knitted and photographed by author, pattern by Norah Gaughan for Berrocco \cite{NG}.}
\end{center}
\end{figure}

In 2014, the knitwear designer Norah Gaughan (known for the intricate mathematics of her cable designs and for ingenious garment constructions) envisaged a hat \textit{Orpheus}, the pattern for which is administered by Berrocco \cite{NG}. The hat is formed from a simple rectangle, but not in the  so-called teabag style with two straight seams. Rather, the rectangle is seamed to itself in an unusual way, as a pita-form, thus defining a 3D space. Of course, if it were seamed all the way around, it could not be worn as a hat! The edges are attached part-way, to create a cavity with one opening, which can then be placed on one's head. The suggested stripes of the knitted fabric, interacting with the single kinked seam, undoubtedly contribute further to the disconcerting effect of this design. A hat made to this pattern is shown in Figure \ref{orpheus}.

\section{Biscornu}
\subsection{Biscornu as textile D-forms}
Biscornu are ornamental cushions, made from two off-set squares of the same size \cite{Bis}. They take their name from a French adjective meaning irregular or more literally multi-horned (c.f. unicorn or tricorn), and are usually made from patch-work fabric, or from even-weave fabric embroidered using a counted-thread technique such as blackwork or cross-stitch, possibly further decorated with beads.  Their popularity with crafters was extremely high for about a decade beginning around 2005. They are frequently used as pincushions, or as hanging decorations. Their origin is obscure, and they do not appear to have been connected to any handcraft tradition. 

\begin{figure}
\begin{center}
{\resizebox*{8cm}{!}{\includegraphics{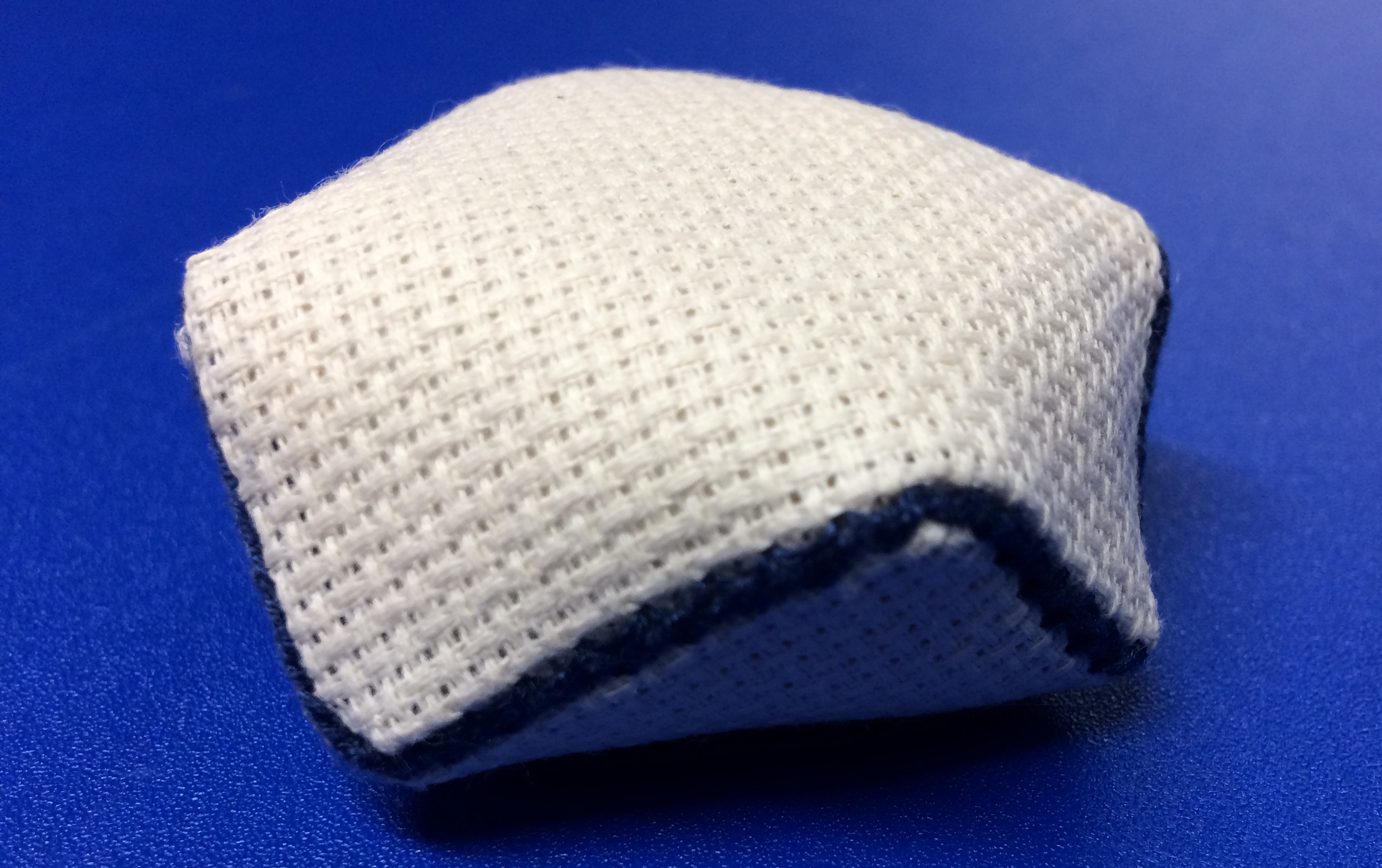}}}
\caption{\label{plain} An undecorated biscornu, made from two squares of Aida cloth, side length 6 cm. }
\end{center}
\end{figure}

A biscornu is a D-form in the original sense of Wills \cite{Wi}, since the boundary of a square is not smooth. The two squares are seamed together with each vertex of one square coinciding with the midpoint of an edge of the other, and \textit{vice versa}.  As can be seen in Figure \ref{plain} the seam does not lie flat, but takes on a zigzag formation, like the crown formation of the molecule {cyclo}-octasulfur. 

As is the case for {cyclo}-octasulphur, the symmetry group of the biscornu  is that of the square anti-prism which it resembles: $D_{4d}$. Smaller square regions that have as their vertices the midpoints of the original edges play the role of the square faces of the anti-prism, and the corners of the original squares bend over, without forming a sharp folded edge, to provide the triangles between the squares (again, see Figure \ref{plain}).  Unsupported, the fabric will not hold this shape, and the cushion is filled in some way. The inherent irregularity is often heightened by pulling the two squares together with a button sewn at the centre of each \cite{Bis}, but this is not essential, and does distort them from being true D-forms. Over-stuffing also needs to be avoided, lest they look more like balls and less like quirky anti-prisms with gently rounded edges. 
\subsection{An axial point group sampler}
In two dimensions, samplers of rosette, frieze and wallpaper groups have been fashioned using various counted thread techniques (for a survey, see \cite{SG}). This paper concludes with images of eleven biscornu designed and constructed to showcase the three-dimensional symmetry group of a biscornu and each of the possible subgroups, including the trivial subgroup; collectively they comprise a single symmetry sampler of this axial point group. The traditional Japanese sashiko stitching form \textit{hitomezashi} \cite{SB,HS} has been chosen as the embroidery form used.

\begin{figure}
\centering
\begin{tikzpicture}
\draw(0,0) rectangle (4,4);
\node at (-0.15,0){$4$};
\node at (4.15,0){$1$};
\node at (-0.15,4){$3$};
\node at (4.15,4){$2$};
\fill[black](2,2) circle (1mm);
\draw (0,2.83)--(-.83,2)--(0,1.17);
\draw(1.17,0)--(2,-0.83)--(2.83,0);
\draw (4,2.83)--(4.83,2)--(4,1.17);
\draw(1.17,4)--(2,4.83)--(2.83,4);
\draw[dashed] (0,2.83)--(1.17,4);
\draw[dashed] (2.83,4)--(4,2.83);
\draw[dashed] (2.83,0)--(4,1.17);
\draw[dashed] (1.17,0)--(0,1.17);
\node at (-1,2){$\gamma$};
\node at (5,2){$\alpha$};
\node at (2,5){$\beta$};
\node at (2,-1){$\delta$};
\node at(2,-2){(a)};
\end{tikzpicture}
\quad
 \begin{tikzpicture}
\draw(0,0) rectangle (4,4);
\node at (-0.15,0){$\gamma$};
\node at (4.15,0){$\delta$};
\node at (-0.15,4){$\beta$};
\node at (4.15,4){$\alpha$};
\fill[black](2,2) circle (1mm);
\draw (0,2.83)--(-.83,2)--(0,1.17);
\draw(1.17,0)--(2,-0.83)--(2.83,0);
\draw (4,2.83)--(4.83,2)--(4,1.17);
\draw(1.17,4)--(2,4.83)--(2.83,4);
\draw[dashed] (0,2.83)--(1.17,4);
\draw[dashed] (2.83,4)--(4,2.83);
\draw[dashed] (2.83,0)--(4,1.17);
\draw[dashed] (1.17,0)--(0,1.17);
\node at (-1,2){$3$};
\node at (5,2){$1$};
\node at (2,5){$2$};
\node at (2,-1){$4$};
\node at(2,-2){(b)};
\end{tikzpicture}
    \caption{(a) The two squares comprising a biscornu are indicated as seen from above before seaming, and the eight vertices are labelled; the square labelled with Greek letters lies under the square whose vertices are numbered. (b) The roto-reflection operation has been applied to (a), the axis of the rotation being perpendicular to the plane of the page at the position shown by the black dot.}
    \label{roto}
\end{figure}
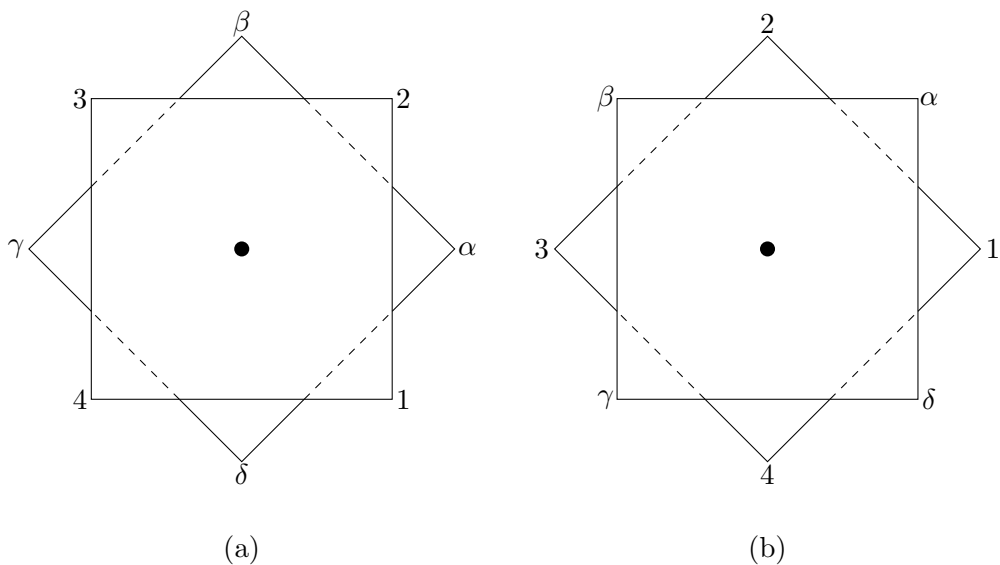

The symmetry group of biscornu, the dihedral group $D_{4d}$, has order 16. One of the group generators is a \textit{rotoreflection} (of order 8): rotation through $\frac{\pi}{4}$ about an axis through the centre of both squares, together with reflection across a plane halfway between these centres and perpendicular to the rotation axis. The effect of this operation is shown schematically in Figure \ref{roto}.

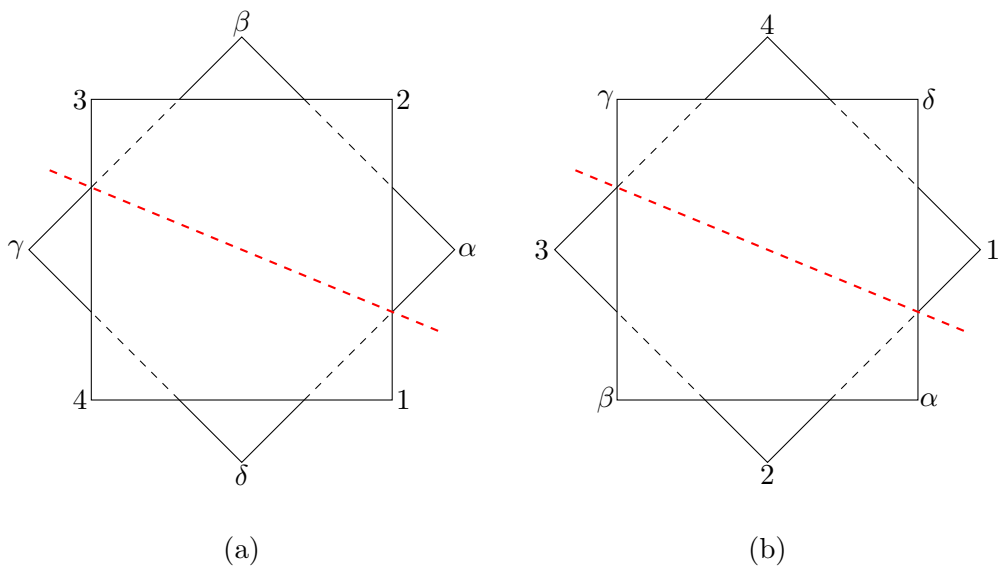
\begin{figure}
    \centering
    \begin{tikzpicture}
\draw(0,0) rectangle (4,4);
\node at (-0.15,0){$4$};
\node at (4.15,0){$1$};
\node at (-0.15,4){$3$};
\node at (4.15,4){$2$};
\draw (0,2.83)--(-.83,2)--(0,1.17);
\draw(1.17,0)--(2,-0.83)--(2.83,0);
\draw (4,2.83)--(4.83,2)--(4,1.17);
\draw(1.17,4)--(2,4.83)--(2.83,4);
\draw[dashed] (0,2.83)--(1.17,4);
\draw[dashed] (2.83,4)--(4,2.83);
\draw[dashed] (2.83,0)--(4,1.17);
\draw[dashed] (1.17,0)--(0,1.17);
\draw[dashed,thick,red] (4.61,0.92)--(-0.61,3.08);
\node at (-1,2){$\gamma$};
\node at (5,2){$\alpha$};
\node at (2,5){$\beta$};
\node at (2,-1){$\delta$};
\node at(2,-2){(a)};
\end{tikzpicture}
\quad
\begin{tikzpicture}
\draw(0,0) rectangle (4,4);
\node at (-0.15,0){$\beta$};
\node at (4.15,0){$\alpha$};
\node at (-0.15,4){$\gamma$};
\node at (4.15,4){$\delta$};
\draw (0,2.83)--(-.83,2)--(0,1.17);
\draw(1.17,0)--(2,-0.83)--(2.83,0);
\draw (4,2.83)--(4.83,2)--(4,1.17);
\draw(1.17,4)--(2,4.83)--(2.83,4);
\draw[dashed] (0,2.83)--(1.17,4);
\draw[dashed] (2.83,4)--(4,2.83);
\draw[dashed] (2.83,0)--(4,1.17);
\draw[dashed] (1.17,0)--(0,1.17);
\draw[dashed,thick,red] (4.61,0.92)--(-0.61,3.08);
\node at (-1,2){$3$};
\node at (5,2){$1$};
\node at (2,5){$4$};
\node at (2,-1){$2$};
\node at(2,-2){(b)};
\end{tikzpicture}
    \caption{(a) A rotation axis generating the other symmetries of the biscornu is indicated by a red dashed line; (b) The result of rotating through $\pi$ about the red axis - note in particular that the chirality of the labelling is changed.}
    \label{diag}
\end{figure}

The other symmetry group generator, of order 2, is a rotation through $\pi$; it is about a line that lies in the reflection plane mentioned above and that makes angle $\frac{\pi}{8}$ with the symmetry axes of the two squares, as shown in Figure \ref{diag} (a). This operation, a non-axial rotation, changes the labelling order of the vertices from anti-clockwise to clockwise, as shown in Figure \ref{diag} (b).

\begin{figure}
\begin{center}
{\resizebox*{8cm}{!}{\includegraphics{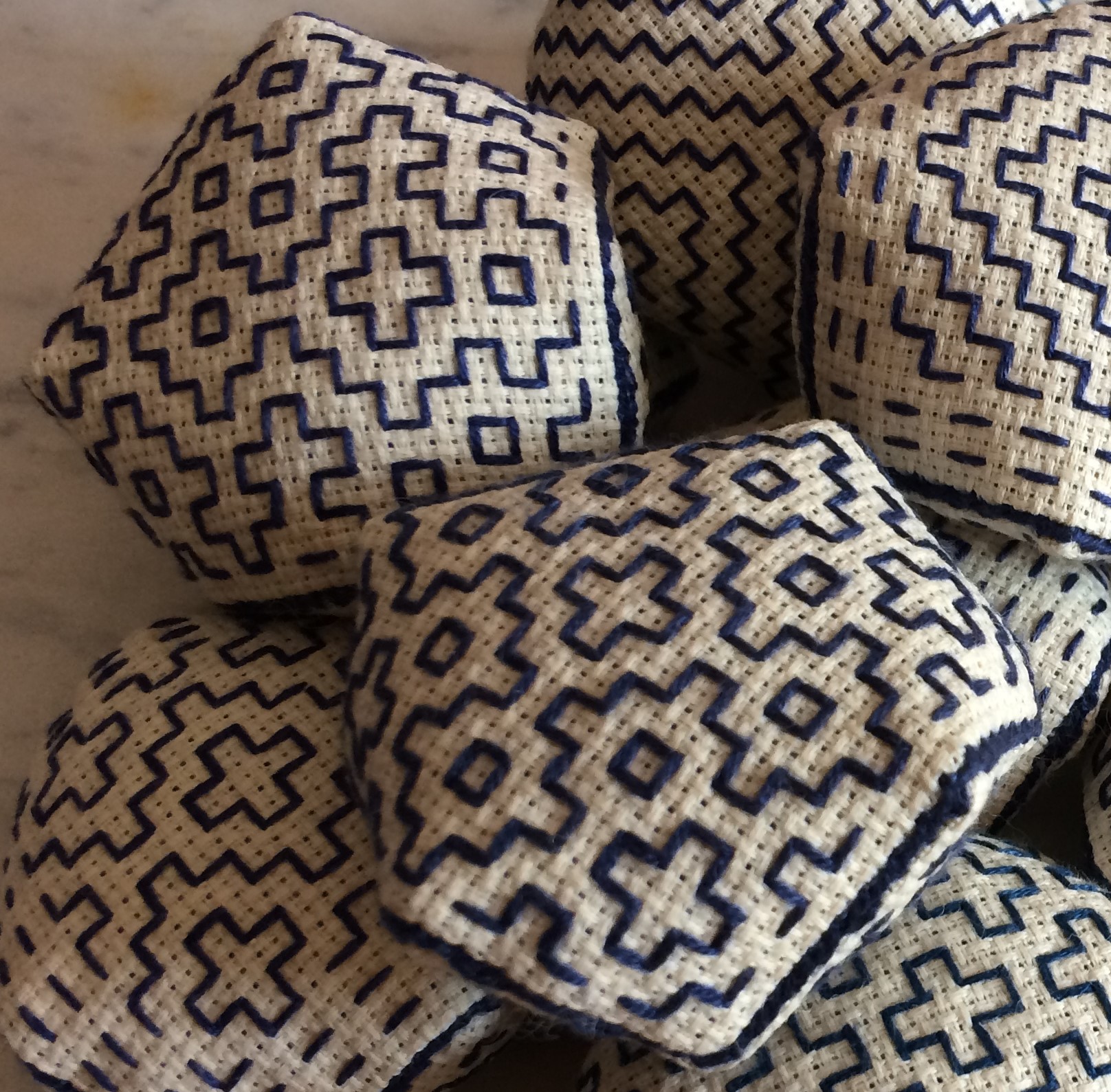}}}
\caption{\label{D4d} In the foreground, a biscornu decorated to retain the full $D_{4d}$ symmetry, with the hidden surface bearing the same motif; behind and to the left, a biscornu `symmetric' only under the identity operator, representing the trivial subgroup. }
\end{center}
\end{figure}

The unadorned biscornu of Figure \ref{plain} has the full $D_{4d}$ symmetry, but so too does the foremost object in  Figure \ref{D4d}. Both squares comprising this biscornu have been decorated with the same Fibonacci snowflake \cite{RRC} which has $D_4$ symmetry. The trivial subgroup has been depicted (the object slightly behind and to the left in Figure \ref{D4d}) by using an aleatoric design; that is, the running stitches have been placed according to the toss of a coin, resulting in no particular symmetry on either square.

There are nine distinct non-trivial subgroups, three each of order 8,  4 and 2.

\begin{figure}
\begin{center}
{\resizebox*{12cm}{!}{\includegraphics{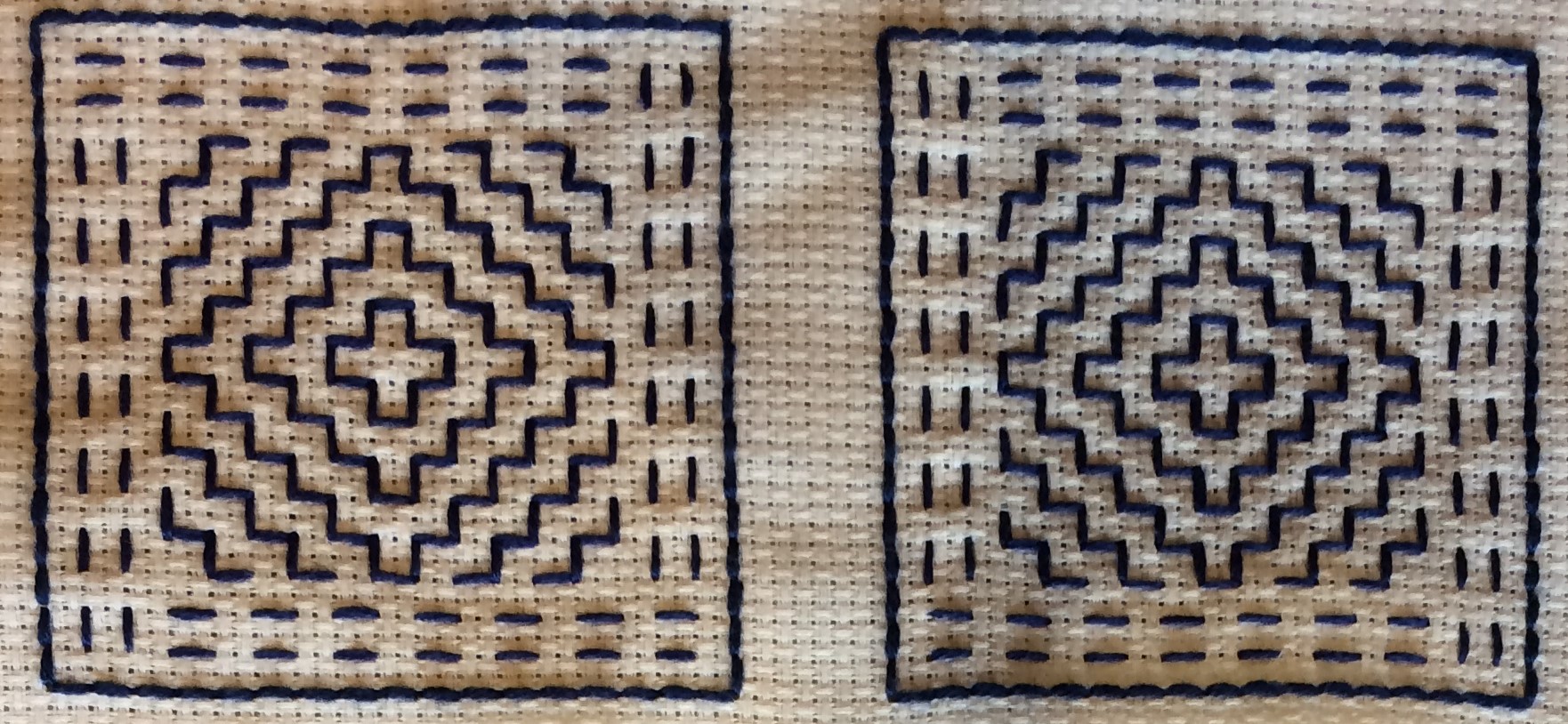}}}
\caption{\label{S8} The two squares which comprise a biscornu symmetric under the order 8 rotoreflection. }
\end{center}
\end{figure}

The first subgroup of order 8 is denoted $S_8$, where in this context the notation refers to eight-fold rotoreflection, not to be confused with the symmetric group on eight elements, or indeed, the chemical symbol for octasulfur. It is isomorphic to $Z_8$ or $C_8$. The two squares comprising the relevant biscornu are shown in Figure \ref{S8}. The centre of the design is the traditional \textit{yamagata} (mountain form). The running stitch borders have \textit{opposite} chirality on the two squares; under the rotoreflection
one maps to the other.

\begin{figure}
\begin{center}
{\resizebox*{12cm}{!}{\includegraphics{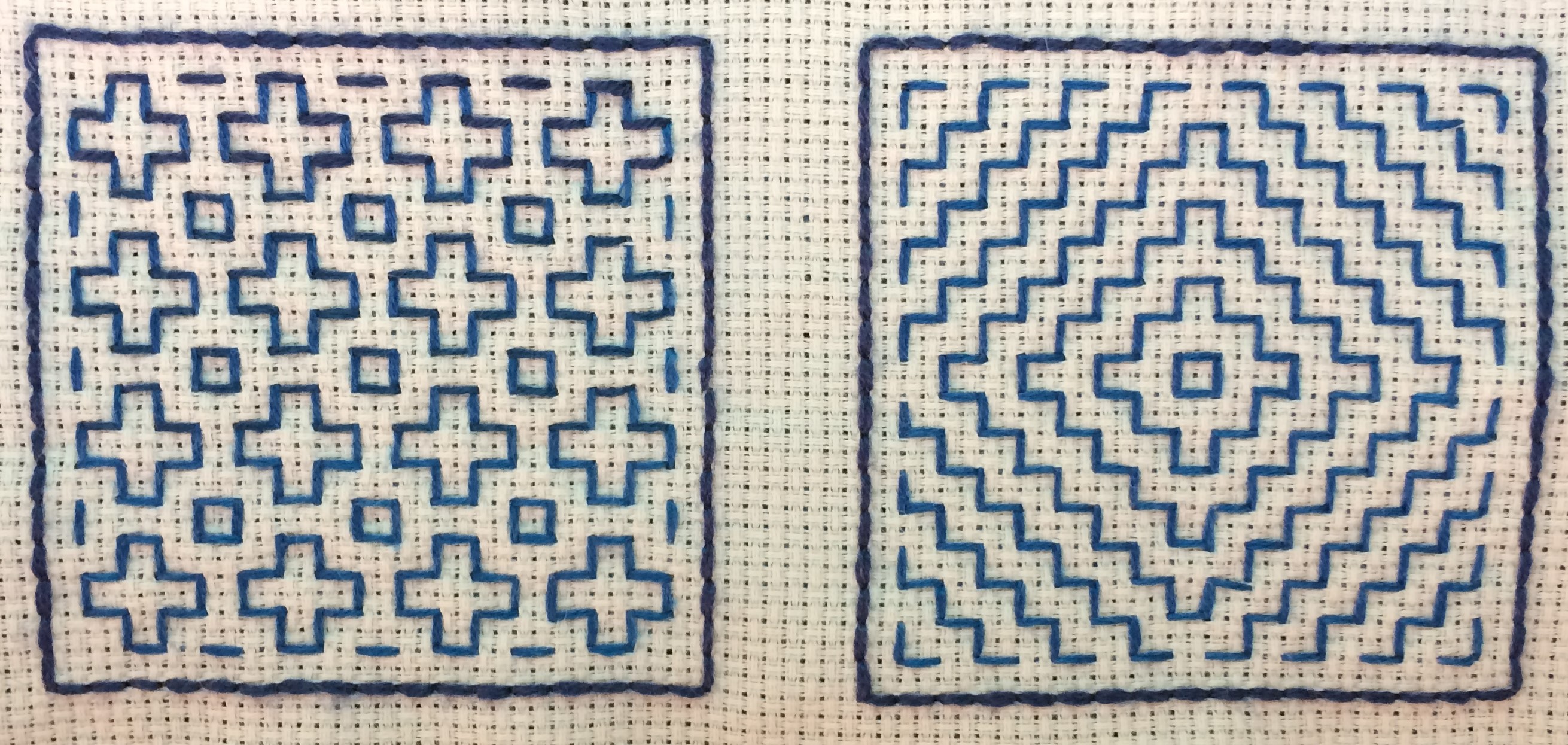}}}
\caption{\label{D4} Each square has $D_4$ symmetry, with four reflection axes, two running with the fabric grid and two diagonally. When the squares are offset and seamed together, the diagonal axes of one line up with the `square' reflection axes of the second. }
\end{center}
\end{figure}

\begin{figure}
\begin{center}
{\resizebox*{8cm}{!}{\includegraphics{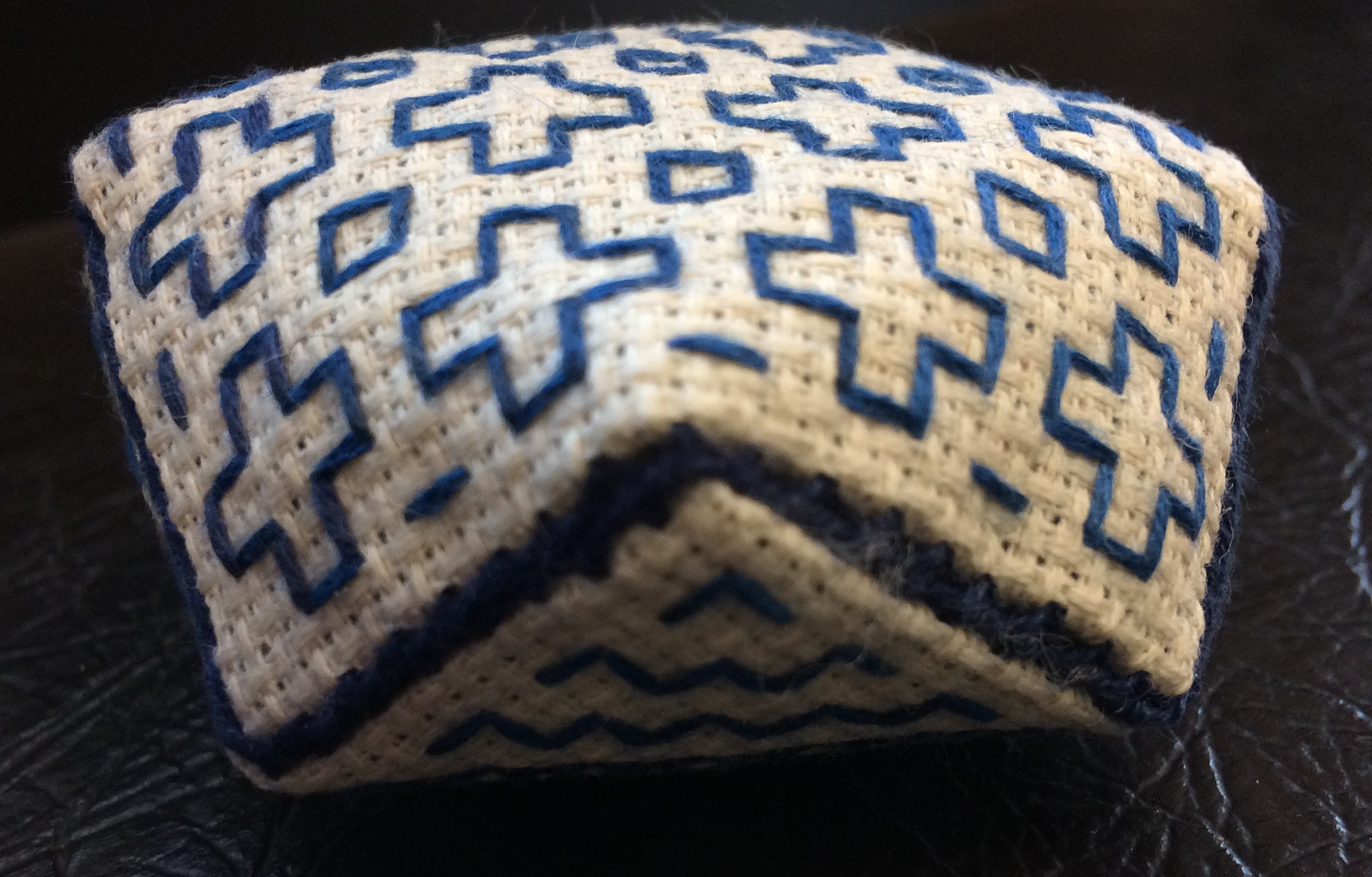}}}
\caption{\label{D4up}  When the two squares from Figure \ref{D4} are seamed together, the overall object has four reflections and four-fold rotation about the axis. This object represents the subgroup $D_4$.}
\end{center}
\end{figure}

The two constituent squares have been decorated with different patterns of the same symmetry to create a biscornu to depict $D_4$ (see Figures \ref{D4} and \ref{D4up}). On one square,  \textit{yamagata} extends over the whole surface. On the other,  the pattern comprises \textit{j\={u}jizashi} (ten-cross, the character for the number ten resembling a cross) and little squares.

\begin{figure}
\begin{center}
{\resizebox*{8cm}{!}{\includegraphics{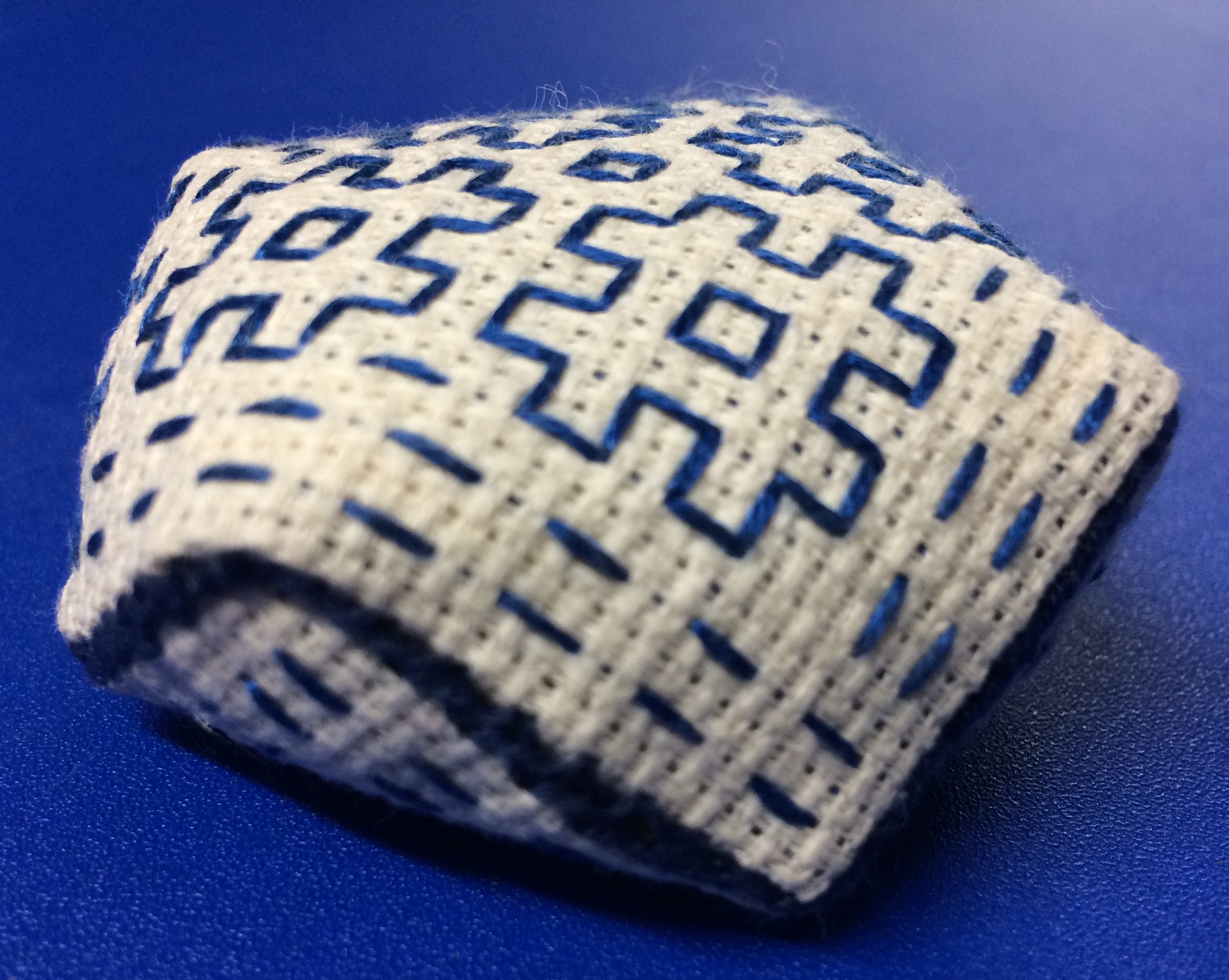}}}
\caption{\label{C4v} A biscornu with $C_{4v}$ symmetry. That the running stitches map correctly when the non-axial rotation of Figure \ref{diag} is applied can be verified in this photo which highlights an edge, not a face. }
\end{center}
\end{figure}

The final order 8 subgroup depicted is $C_{4v}$, shown in Figure \ref{C4v}. The well-kerb (\textit{igetazashi}) design has been surrounded by running stitches of the \textit{same} chirality on each square. This object is unchanged by rotation about the central axis through $\frac{\pi}{2}$ (order 4)  and by the non-axial rotation (order 2). 

\begin{figure}
\begin{center}
{\resizebox*{12cm}{!}{\includegraphics{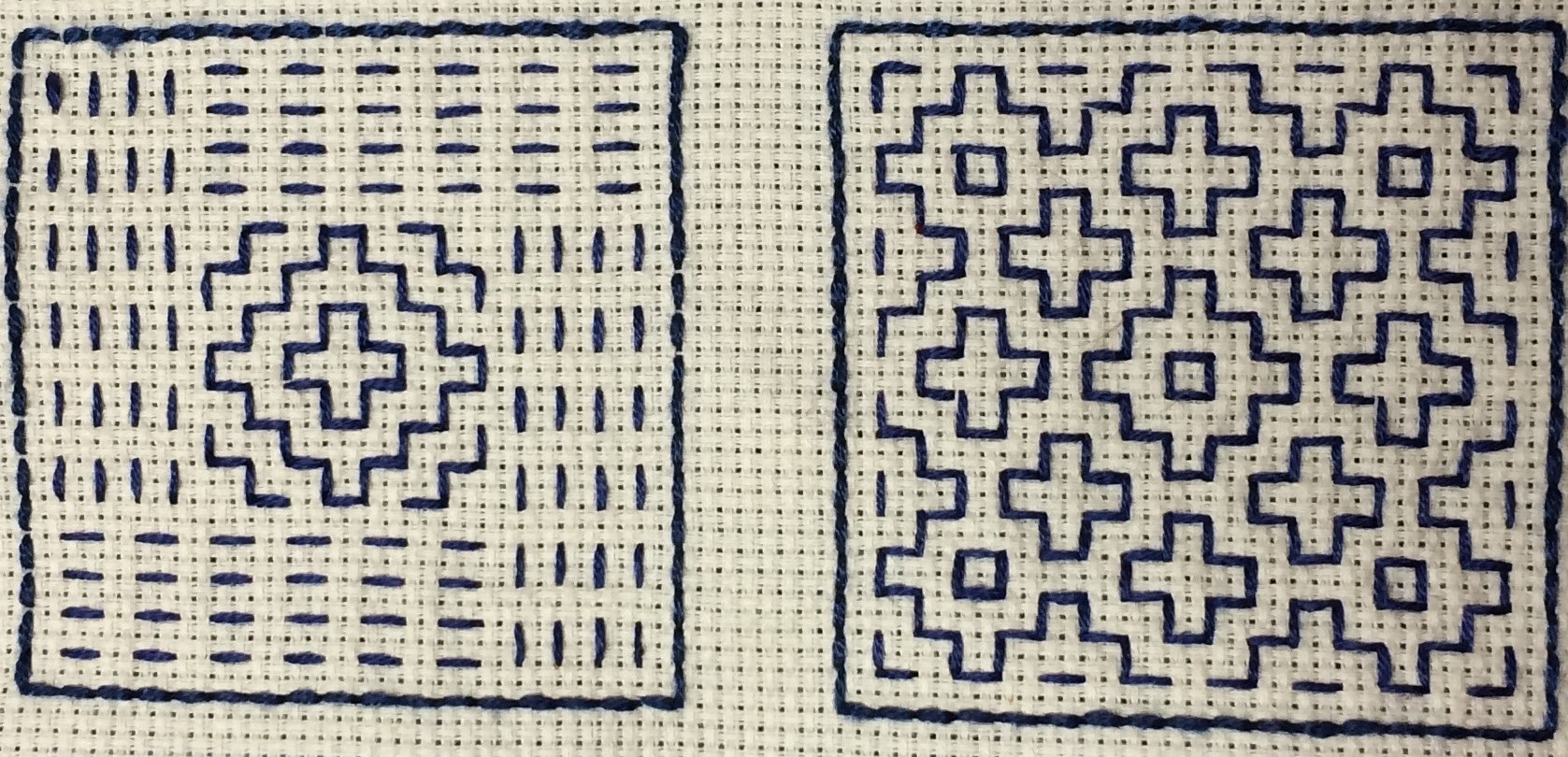}}}
\caption{\label{C4flat} A square with the desired symmetry (left) is to be combined with a square that has additional symmetry, to create a biscornu to illustrate $C_4$. }
\end{center}
\end{figure}

To illustrate the first of the order 4 subgroups $C_4$, as shown in Figure \ref{C4flat} one square bears a design with that symmetry (a variant of \textit{kakinohanazashi} persimmon flower stitch, surrounded by wide chiral bands of running stitch) while the other has further symmetry ($D_4$). This second square features a pleasing design of squares, crosses and a persimmon flower; in fact it is essentially the reverse fabric (or dual \cite{HS}) of the Fibonacci snowflake in Figure \ref{D4d}. 

\begin{figure}
\begin{center}{\resizebox*{12cm}{!}
{\includegraphics{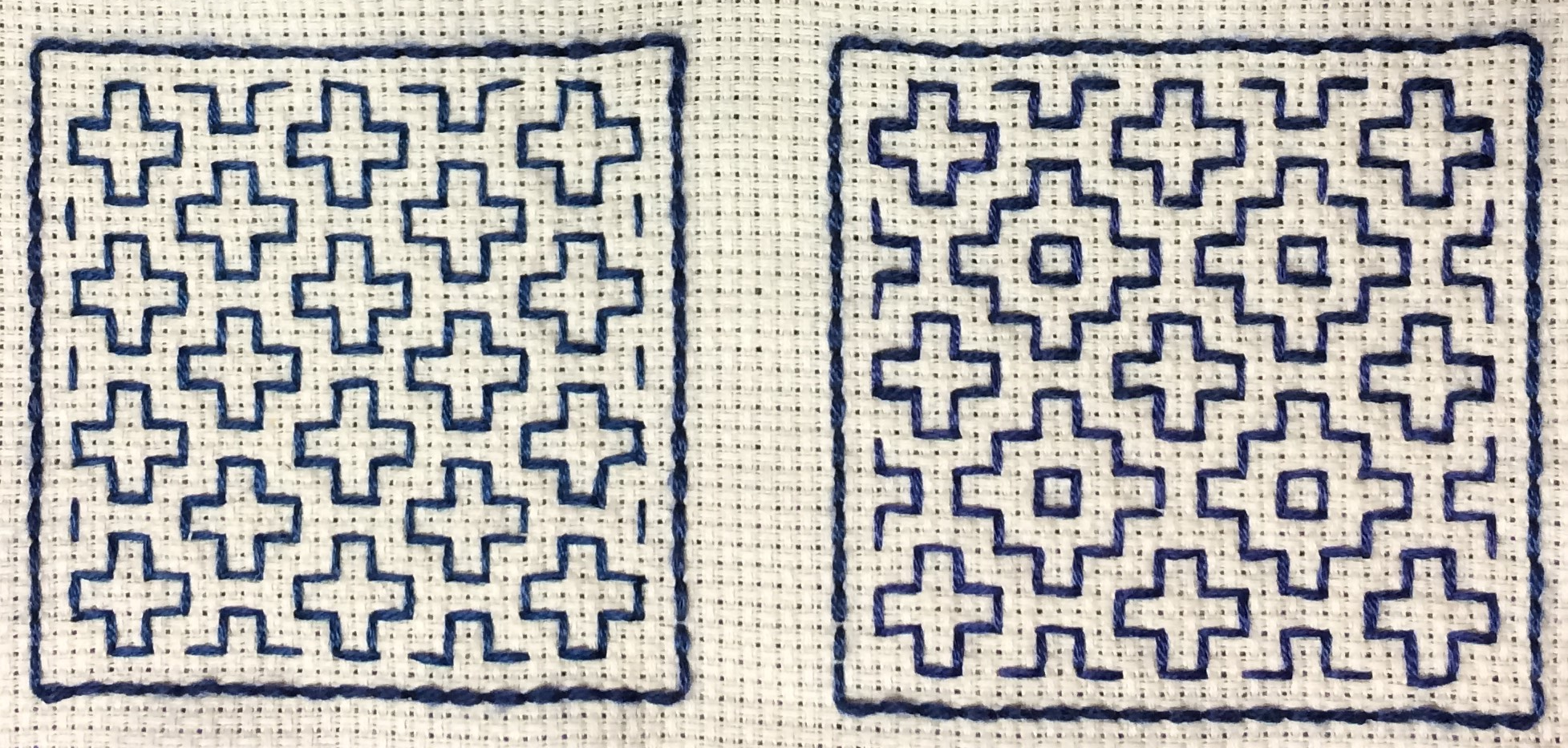}}}
\caption{\label{D2} In this case, combining a square with the desired symmetry (left) with a square that has additional symmetry, a biscornu to illustrate $D_2$ is created. }
\end{center}
\end{figure}

A similar approach has been used to design a biscornu that has $D_2$ symmetry. As shown in Figure \ref{D2} a close-packed arrangement of \textit{j\={u}jizashi}  has been paired with an square that has additional symmetry in its design. When this second square is offset, its two diagonal reflection axes line up with those of the other square, and being symmetric under rotation through $\frac{\pi}{4}$, it is also symmetric under the rotation through $\frac{\pi}{2}$.

\begin{figure}
\begin{center}
{\resizebox*{12cm}{!}{\includegraphics{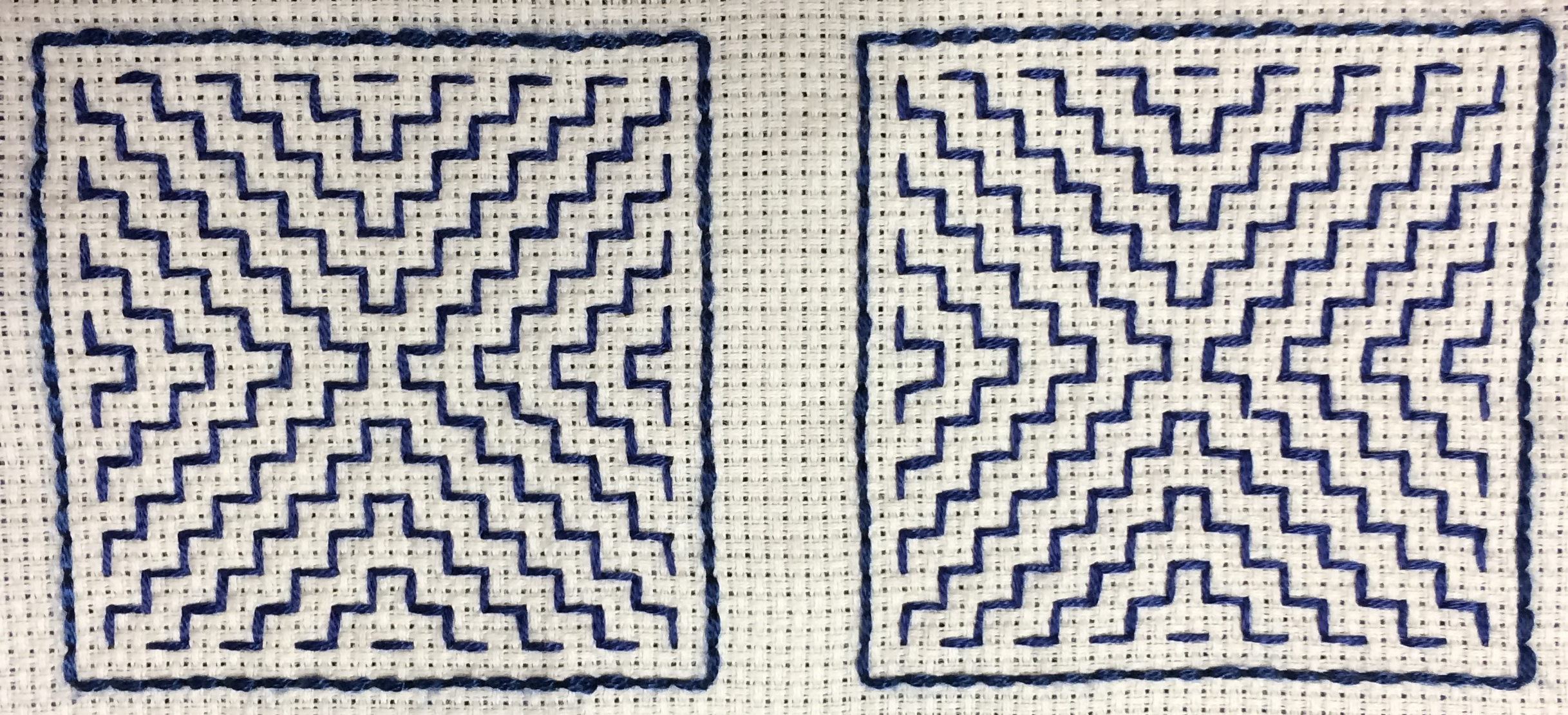}}}
\caption{\label{C2v} Two identical designs which use the mountain form stitch and which combine to give a biscornu with $C_{2v}$ symmetry.}
\end{center}
\end{figure}

The final subgroup of order four $C_{2v}$ is generated by rotation through $\pi$ about the principal axis and the non-axial rotation, both of order 2. The squares used to create a biscornu with this symmetry bear identical designs (with $D_2$ symmetry) as shown in Figure \ref{C2v}. The designs do not have diagonal symmetry axes, so that when one square is offset the resultant biscornu would have only $C_2$ symmetry were it not that the squares map to each other under the non-axial rotation. 

\begin{figure}
\begin{center}
{\resizebox*{8cm}{!}{\includegraphics{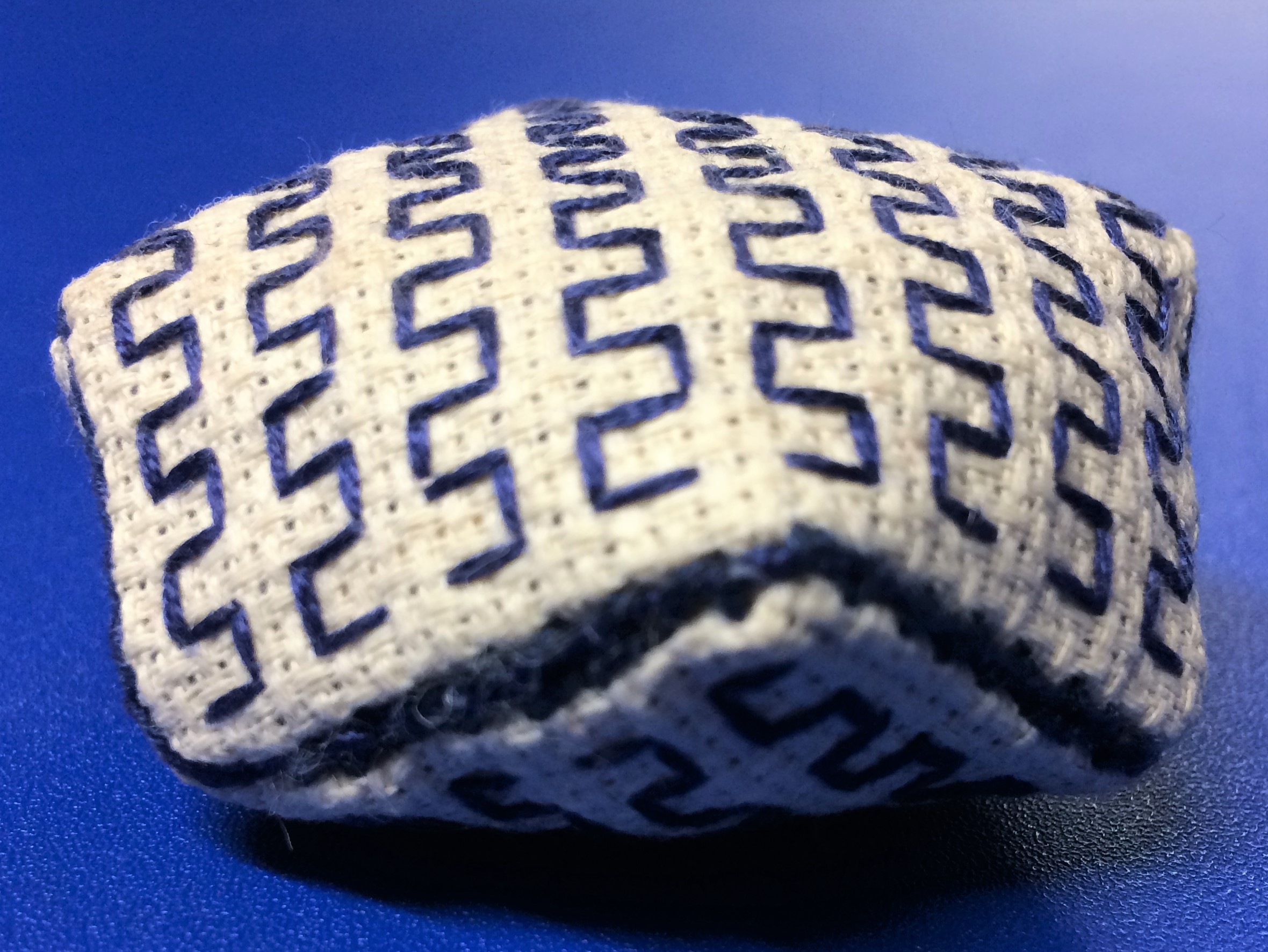}}}
\caption{\label{OC2}  Since the identical squares do not have diagonal symmetry, the resultant biscornu has  the non-axial rotation as its only symmetry. The stitch design is called \textit{hirayama michi} (mountain passes).}
\end{center}
\end{figure}

\begin{figure}
\begin{center}
{\resizebox*{12cm}{!}{\includegraphics{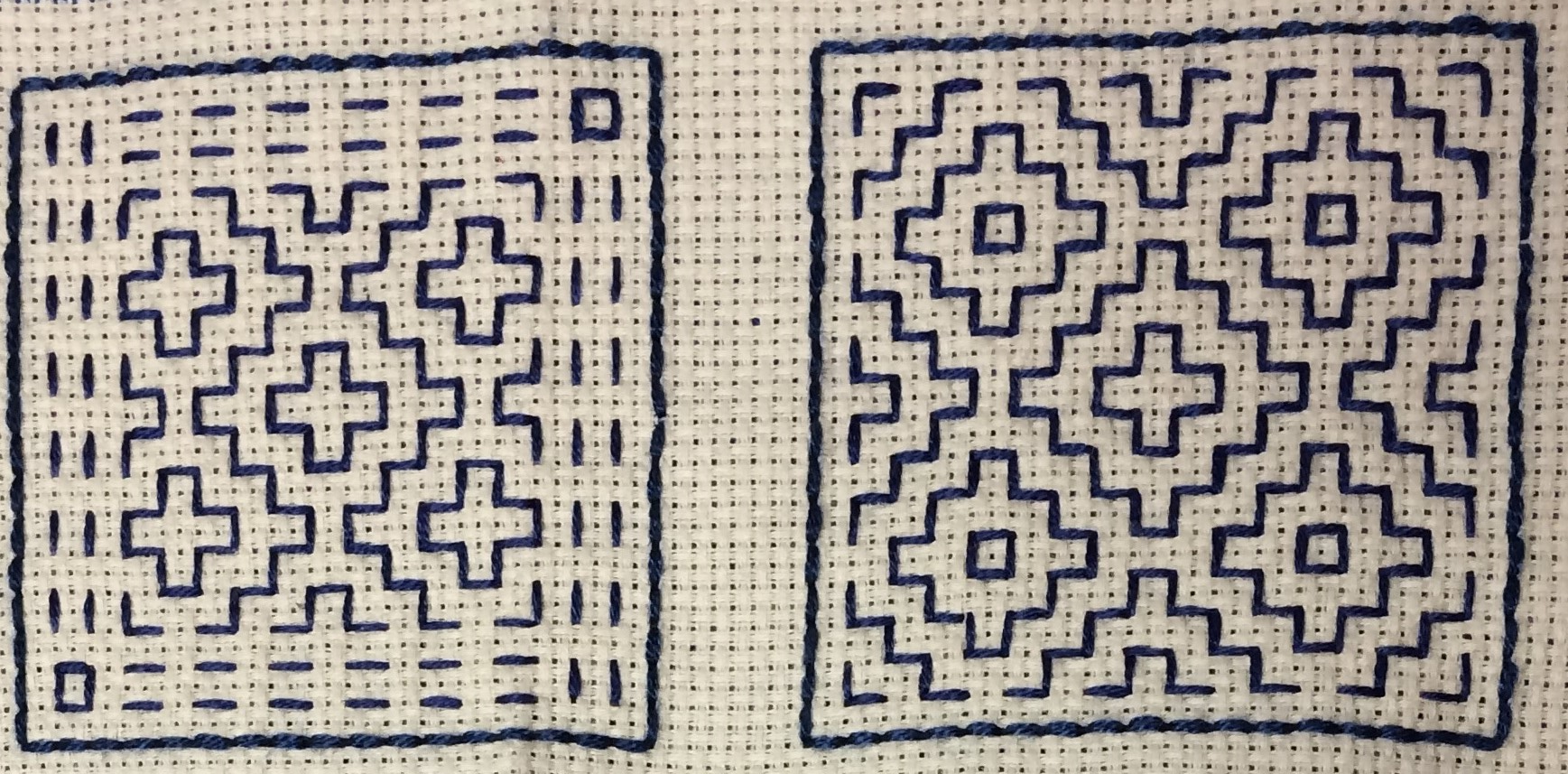}}}
\caption{\label{C2} Both the running stitch frame and the central design in the square on the left are symmetric only under rotation through $\pi$. The by-now familiar approach of using a second square with more symmetry, in this case a design of two kinds of persimmon, has been used to yield a biscornu with bilateral symmetry. }
\end{center}
\end{figure}

\begin{figure}
\begin{center}
{\resizebox*{12cm}{!}{\includegraphics{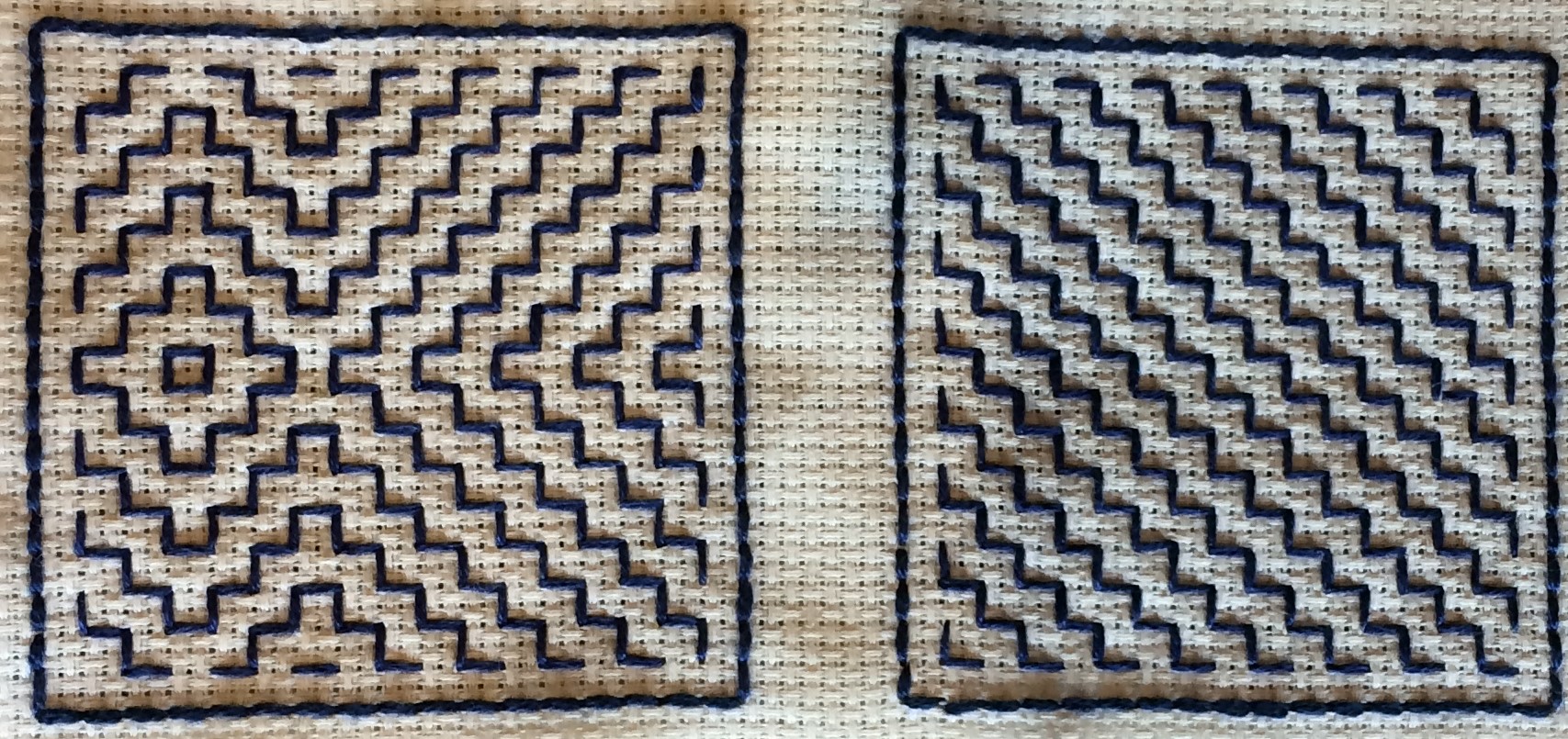}}}
\caption{\label{Cs} The \textit{yamagata} form in the square on the left could have been extended; instead, a persimmon has been stitched like a sun above the mountain. The second square can be joined to the first in only one of two ways (not four) in order to yield a $C_s$ biscornu.}
\end{center}
\end{figure}

There are two $C_2$ subgroups of $D_{4d}$. One is generated by the non-axial rotation (only) and the other by rotation through $\pi$ about the principal axis. The first of these is shown in Figure \ref{OC2}, and the constituent squares of the second in Figure \ref{C2}.
The final subgroup, also of order 2,  is $C_s$ and corresponds to an object with a single bilateral symmetry. To create such a biscornu, a square with a step design \textit{dan tsunagi}, that has a single diagonal symmetry axis, has been combined with a design having a single symmetry axis parallel to one pair of sides. The squares are shown in Figure \ref{Cs}.

\section{Conclusion}
While intended to be artistically interesting in their own right, these decorated D-forms with their symmetric designs could be useful as manipulables in the teaching of algebra or of molecular chemistry. I can attest that devising stitch patterns that have the desired symmetries  certainly improves one's understanding of three dimensional geometry, as well as being creatively satisfying.

\section*{Acknowledgements} With thanks to Bruce Atkinson for permitting me to use photos of his designs, and for an entertaining correspondence.

\end{document}